\documentclass{article}

     \usepackage[final]{ml4eng_2020}
\usepackage[utf8]{inputenc} 
\usepackage[T1]{fontenc}    
\usepackage{hyperref}       
\usepackage{url}            
\usepackage{booktabs}       
\usepackage{amsfonts}       
\usepackage{nicefrac}       
\usepackage{microtype}      
\usepackage{graphicx} 
\usepackage{algorithm}
\usepackage[noend]{algpseudocode}
\usepackage{mathtools}

\title{On Training Effective Reinforcement Learning Agents for Real-time Power Grid Operation and Control}

\author{%
  Ruisheng Diao\thanks{Dr. Ruisheng Diao is with GEIRI North America as Deputy Department Head, AI \& System Analytics, leading the development of a number of AI-based applications in power systems (autonomous voltage control, line flow control, intelligent maintenance scheduling for power utilities and control centers). Cell phone: (480)-414-7095, website: https://www.linkedin.com/in/ruisheng-diao-ph-d-pe-789a9655/} \\
  AI \& System Analytics\\
  GEIRI North America\\
  San Jose, CA 95134 \\
  \texttt{ruisheng.diao@geirina.net} \\
  \And
   Di Shi \\
   AI \& System Analytics\\   
   GEIRI North America \\
   San Jose, CA \\
   \texttt{di.shi@geirina.net} \\
   \And
   Bei Zhang \\
   AI \& System Analytics\\   
   GEIRI North America \\
   San Jose, CA \\
   \texttt{beizhangzl@gmail.com} \\
   \And
   Siqi Wang \\
   AI \& System Analytics\\   
   GEIRI North America \\
   San Jose, CA \\
   \texttt{siqi.wang@gmail.com} \\   
   \And
   Haifeng Li \\
   System Operations\\   
   Jiangsu Electric Power Company \\
   Nanjing, China \\
   \texttt{lhfcloud@163.com} \\   
   \And
   Chunlei Xu \\
   Power System Automation\\   
   Jiangsu Electric Power Company \\
   Nanjing, China \\
   \texttt{chunleixu@js.sgcc.com.cn} \\     
   \And
   Tu Lan \\
   AI \& System Analytics\\  
   GEIRI North America \\
   San Jose, CA \\
   \texttt{lantusky2017@gmail.com} \\    
   \And
   Desong Bian \\
   AI \& System Analytics\\  
   GEIRI North America \\
   San Jose, CA \\
   \texttt{desong.bian@geirina.net} \\     
   \And
   Jiajun Duan \\
   AI \& System Analytics\\  
   GEIRI North America \\
   San Jose, CA \\
   \texttt{jiajun.duan@geirina.net} \\     
 
}

\begin{document}

\maketitle
\begin{abstract}
Deriving fast and effectively coordinated control actions remains a grand challenge affecting the secure and economic operation of today’s large-scale power grid. This paper presents a novel artificial intelligence (AI) based methodology to achieve multi-objective real-time power grid control for real-world implementation. State-of-the-art off-policy reinforcement learning (RL) algorithm, soft actor-critic (SAC) is adopted to train AI agents with multi-thread offline training and periodic online training for regulating voltages and transmission losses without violating thermal constraints of lines. A software prototype was developed and deployed in the control center of SGCC Jiangsu Electric Power Company that interacts with their Energy Management System (EMS) every 5 minutes. Massive numerical studies using actual power grid snapshots in the real-time environment verify the effectiveness of the proposed approach. Well-trained SAC agents can learn to provide effective and subsecond (<20 ms) control actions in regulating  voltage profiles and reducing transmission losses.
\end{abstract}

\section{Introduction}
Over recent years, the electricity sector has undergone significant changes, with ever-increasing penetration of intermittent energy resources, storage devices, and power electronics equipment integrated into decades-old grid infrastructure, causing more stochastic and dynamic behavior that affects the secure and economic operation of the grid. Various control measures exist to restore the imbalance of active and reactive power to ensure voltage, frequency and line flows operating within their normal ranges. Compared to mature controllers that function well at local levels, deriving system-wide optimal controls by coordinating many controllers for real-time operating conditions while complying with all security constraints remains a grand challenge. 

Such a control problem is known as security constrained AC optimal power flow (AC OPF), where control objectives can be minimization of total generation costs, total network losses, amount of control actions or a combination of all while respecting physics and security constraints. Many prior research efforts were proposed to solve the large-scale non-convex AC OPF problem considering various constraints with mixed-integer variables, including nonlinear programming, quadratic programming, Lagrangian relaxation, interior point method and mixed-integer linear programming. A comprehensive survey is conducted in [1] that summarizes the main challenges and state-of-the-art techniques. In [2]-[3], semidefinite programming (SDP) relaxation-based methods are proposed for solving the multi-objective OPF problem. In [4], the authors proposed a quasi-Newton method using second-order information for solving real-time OPF problem. However, due to the non-convex and complex nature of the problem (NP-Hard) and the need for accurate system-wide models in real-time environment, deriving solutions for the AC OPF problem with all security constraints is very challenging. Thus, relaxation of constraints and simplifications of system models, e.g., DC based OPF, are typically used to obtain feasible and suboptimal solutions. In fact, nearly all today's vendors' tools adopt DC-based models for obtaining fast OPF solutions in power industry. To effectively deal with the new challenges and derive fast OPF solutions, new approaches are much-needed [5].

Recent success of applying RL technologies in various control problems including game of GO, self-driving cars and robotics provides a promising way of deriving instant control actions for reaching feasible and suboptimal AC OPF solutions. Recently, RL-based algorithms were reported for regulating voltage profiles [6] and transmission line flows [7], load-frequency control [8], and emergency control for enhancing transient voltage behavior [9]. In this paper, a novel method is presented to perform multi-objective power grid control in real time, by adopting the state-of-the-art off-policy RL algorithm, soft actor-critic, with superior performance in convergence and robustness over other RL algorithms [10]. The novelty of this method includes: (1) the AC OPF control problem is formulated as Markov Decision Process (MDP) where RL-based algorithms can be applied to derive suboptimal solutions; (2) it provides a general and flexible framework to include various types of control objectives and constraints of a power grid when training AI agents; (3) once properly trained, RL agents can provide subsecond control actions to regulate bus voltages and transmission losses when sensing abnormal power system states, in case of abrupt changes in voltages and line flows; and (4) multi-threading-based training process of SAC agents with periodical model update is developed to ensure long-term control effectiveness and mitigate the over-fitting issue. The remainder of this paper is organized as follows. Section 2 provides the proposed RL-based methodology for deriving real-time control actions including architecture design, definition of RL components, and implementation details in a control center. In Section 3, case studies are conducted on a city-level power network for demonstrating the effectiveness of the proposed approach. Finally, conclusions are drawn in Section 4 with future work identified.

\section{Proposed Solution for Real-time Multi-objective Power Grid Control Using Soft Actor-Critic}

\subsection{Problem Formulation}
Without losing generality, this work mainly targets at deriving real-time corrective operational control decisions for power grid operating conditions at an interval of 5 minutes in a control center. The control objectives include regulating bus voltages within secure zones and minimizing transmission losses while respecting power flow equations and physical constraints, e.g., line ratings and limits of generators. The mathematical formulation of the control problem is given below:
\paragraph{Objective function}:
\begin{equation}
\text{minimize} \sum_{i,j}P_{loss}(i,j),  (i,j) \in \Omega_L  \label{eq_1}
\end{equation}
\paragraph{Subject to}:
\begin{equation} \label{eq_2}
\sum\limits_{n\in G_i}P_n^g - \sum\limits_{n\in D_i}P_m^d - g_iV_i^2 =  \sum\limits_{j\in B_i}P_{ij}(y) 
\end{equation}
\begin{equation} \label{eq_3}
\sum\limits_{n\in G_i}Q_n^g - \sum\limits_{n\in D_i}Q_m^d - b_iV_i^2 = \sum\limits_{j\in B_i}Q_{ij}(y) \\
\end{equation}
\begin{equation} \label{eq_4}
P_{ij}(y)=g_{ij}V_i^2-V_iV_j(g_{ij}\cos(\theta_{i}-\theta_{j}) +b_{ij}\sin(\theta_{i}-\theta_{j})), (i,j)\in \Omega_{L} \\
\end{equation}
\begin{equation} \label{eq_5}
Q_{ij}(y)=-V_i^2(b_{ij0}+b_{ij})-V_iV_j(g_{ij}\sin(\theta_{i}-\theta_j) - b_{ij}\cos(\theta_{i}-\theta_j)), (i,j)\in \Omega_{L}
\end{equation}

\begin{equation} \label{eq_6}
P_n^{\min} \leq P_n \leq P_n^{\max},  n \in G;  
\end{equation} 
\begin{equation} \label{eq_7}
Q_n^{\min} \leq Q_n \leq Q_n^{\max},  n \in G 
\end{equation}  

\begin{equation} 
V_i^{\min} \leq V_i \leq V_i^{\max},  i \in B \label{eq_8}
\end{equation}

\begin{equation}
\sqrt{P_{ij}^2 + Q_{ij}^2} \leq S_{ij}^{\max}, (i,j) \in \Omega_L \label{eq_9}
\end{equation} 

where $P_{loss}$ is power losses on transmission line connecting bus (or node) $i$ and bus $j$; $P_{n}^g$ is active power injection into bus $n$; $P_{m}^d$ is active power consumption at bus $m$; $\theta_i$ and $V_i$ are voltage phage angle and magnitude at bus $i$, in polar form, respectively. $b_{ij}$ and $g_{ij}$ are conductance and susceptance of transmission line. $P_{ij}$ and $Q_{ij}$ stand for active and reactive power on a transmission line. 
Eq.~\ref{eq_2} through Eq.~\ref{eq_5} represent  power flow equations, representing quasi-steady-state conditions of a power grid, which need to be met all the time. Eq.~\ref{eq_6} and Eq.~\ref{eq_7} are active and reactive power output limits of each generator, respectively. Eq.~\ref{eq_8} and Eq.~\ref{eq_9} specify bus voltage secure zones and line flow limits to be controlled, respectively.

\subsection{Overall Flowchart of the Proposed Methodology}

\begin{figure}
  \centering
  \includegraphics[width=14.2cm,height=15.2cm]{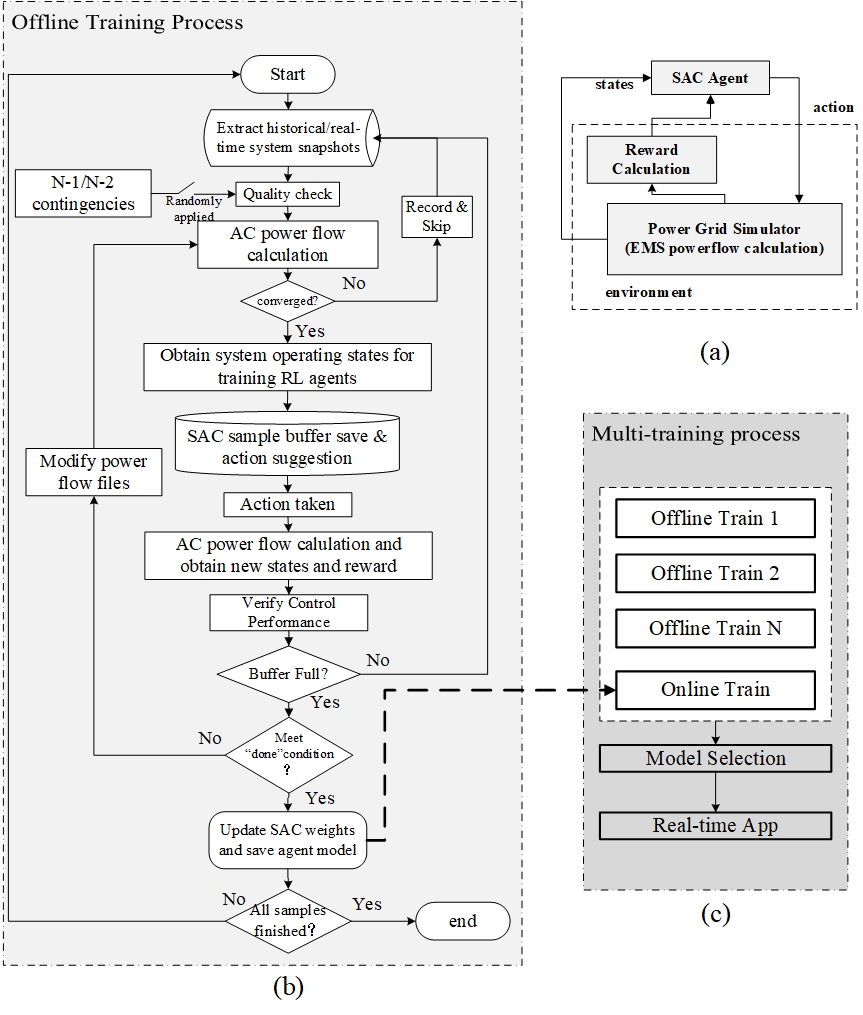}
  \caption{Main flowchart and key function modules of the proposed methodology}
  \vspace{-10pt}
  \label{fig:mainflowchart} 
\end{figure}

Deriving multi-objective real-time control actions respecting both equality and inequality constraints shown in Eq.~\ref{eq_1} through Eq.~\ref{eq_9} can be formulated as a discrete-time stochastic control process, a.k.a., MDP. Among various RL techniques, the off-policy, SAC method is adopted because of its superior performance in fast convergence and robustness, which maximizes the expected reward by exploring as many control actions as possible, leading to a better chance of finding optimum. The detailed algorithm description of SAC can be found in [10]. The main flowchart of the proposed methodology is depicted in Fig~\ref{fig:mainflowchart}, including three key modules:

Module (a): Fig~\ref{fig:mainflowchart}(a) provides the interaction process between the power grid simulation environment (an AC power flow solver satisfying Eq.~\ref{eq_2} through Eq.~\ref{eq_7}) and the SAC agent. The environment receives control actions, outputs the corresponding next system states and calculates the reward; while the SAC agent receives states and reward before outputting control actions, in order to satisfy Eq.~\ref{eq_8}, Eq.~\ref{eq_9} and minimize the objective function, Eq.~\ref{eq_1}.

Module (b): Fig~\ref{fig:mainflowchart}(b) shows the offline training process of an SAC agent. Representative power grid operating snapshots are collected from EMS for preprocessing. System state variables are extracted from those converged snapshots and fed into SAC agent training module, where neural networks are used to establish direct mappings between system states and control actions. The controls are then verified by another AC power flow solution to calculate reward values before updating SAC agent weights for maximizing long-term expected reward and entropy.

Module (c): To ensure long-term effectiveness and robustness of SAC agent, multiple training processes with different sets of hyperparameters are launched simultaneously, including several offline training processes and one online training process (initialized by the best offline-trained model), shown in Fig~\ref{fig:mainflowchart}(c). The best-performing model selected from these processes is then used for application in real-time environment. 
\subsection{Training Effective SAC Agents}
Several key elements in training effective RL agents for multi-objective real-time power flow control are given below:
\subsubsection{Episode, Terminating Conditions, Simulation Environment}
Each episode is defined as a quasi-steady-state operating snapshot, obtained from the EMS system at an interval of 5 minutes. Termination condition of training episodes includes: i) no more voltage or thermal violations with reduction of transmission losses reaching a threshold, e.g., 0.5\%; ii) power flow calculation diverges; or iii) the maximum number of control iteration is reached.
The power grid simulation environment used in this work is the AC power flow calculation program used by system operators in a control center, which typically uses the Newton-Raphson method for obtaining converged solution of system snapshots [11].
\subsubsection{State Space}
The state space is formed by including bus voltage magnitudes ($V$), phase angles ($\theta$), active and reactive power on transmission lines ($P_{ij}$ and $Q_{ij}$), collected from real-time power grid snapshots in the areas to be monitored and controlled. The batch normalization technique is applied to different types of variables for maintaining consistency and improving agent training efficiency. 
\subsubsection{Reward Definition}
The reward value at each control iteration when training SAC agent consider adopts the following logic:
\begin{equation}
reward =  
\begin{cases*}
-(\Delta P_{overflow}) / 10 - (\Delta V_{violation}) / 100, & if voltage or flow violation \\
50- \Delta P_{loss} * 1000, &if $\Delta P_{loss} < 0$ \\
-100, & if $\Delta P_{loss} \geq 2\%$ \\
-1 - \Delta P_{loss} * 50, &otherwise
\end{cases*}
\end{equation}
\begin{equation}
\Delta P_{loss} = P_{loss} - P_{loss\_pre}
\end{equation}
\begin{equation} 
\Delta P_{overflow} = \sum\limits_{i,j}^N((S_{ij}-S_{ij}^{\max})^2 
\end{equation} 
\begin{equation} 
\Delta V_{violation} =\sum\limits_i^M((V_i-V_i^{\max})*(V_i-V_i^{\min})) 
\end{equation} 

where $N$ is the total number of lines with thermal violation; $M$ is the total number of buses with voltage violation; $P_{loss}$ is the present transmission loss value and $P_{loss\_pre}$ is the line loss at the base case. 
\subsubsection{Action Space}
In this work, conventional generators are used to regulate voltage profiles and transmission losses. A control vector is then created to include voltage set points ($V_{set}$) at each selected power plant as continuous values, e.g., [0.9,1.1] p.u. The same control command is applied to all available generators inside the same power plant. 
The details of training SAC agents are given in Algorithm~\ref{sac}, where the gradients of different networks are calculated using equations provided in [10],[12].

\begin{algorithm}
\caption{SAC Training Algorithm for Multi-objective Real-time Power Grid Control} \label{sac}
\begin{algorithmic}
\State Initialize: policy network $\phi$ $\leftarrow$ random weights; \textit{Q} networks $\theta_1$, $\theta_2$ $\leftarrow$ random weights; target networks $\psi_1$ $\leftarrow$ $\theta_1$, $\psi_2$ $\leftarrow$ $\theta_2$; replay buffer \textit{D} $\leftarrow$ $\{\}$
\State Set up power grid simulation environment $env$ 
\For{$episode$ in range ($n\_episodes$)}
\State reset environment and get initial conditions, $state$, $reward$, $P_{loss\_ini}$, $V_{set\_ini}$, $done$ $\leftarrow$ env.reset()
\For{step $t$ in $itertools.count()$}
        \If {replay buffer \textit{D} not full}
            \State random action $a_{t}$ $\leftarrow$ random(A)
        \Else 
            \State action a$_t$ $\leftarrow$ sample from policy $\pi$($a_{t}$\big|$s_{t}$)
        \EndIf

\State  SAC agent interacts with environment, s$_{t+1}$, $r_t$, $done$ $\leftarrow$ env.step(a$_t$,s$_t$)

\State Store transition ($s_{t}$, $s_{t+1}$, $r$, $a_t$, $done$) to $D$

        \If {$done$ == True or $t$ == maxlen}
            \State break
        \EndIf

        \If {$episode$ $\geq$ $batch\_size$}
            \State $\{$($s_{t}$, $a_{t}$, $r$, $s_{t+1}$)$\}$ $\leftarrow$ sample from \textit{D}
            
            \State $\theta_{i}$ $\leftarrow$ $\theta_{i}$ - $\lambda_{Q}$$\nabla_{\theta}$$J_{Q}$($\theta_{i}$), i $\in$ (1,2)

            \State $\phi$ $\leftarrow$ $\phi$ - $\lambda_{\pi}$$\nabla_{\phi}$$J_{\pi}$($\phi$)

            \State $\alpha$ $\leftarrow$ $\alpha$ - $\lambda$$\nabla_{\alpha}$$J(\alpha)$

            \State $\psi_{i}$ $\leftarrow$ $\tau$$\theta_{i}$ + (1 - $\tau$)$\psi_{i}$, i $\in$ (1,2)
        \EndIf
    \EndFor
    \State Save agent model and write logs
\EndFor
\end{algorithmic}
\end{algorithm}

\section{Control Performance of SAC Agents}
\begin{figure}
  \centering
  \includegraphics[width=13.7cm,height=6cm]{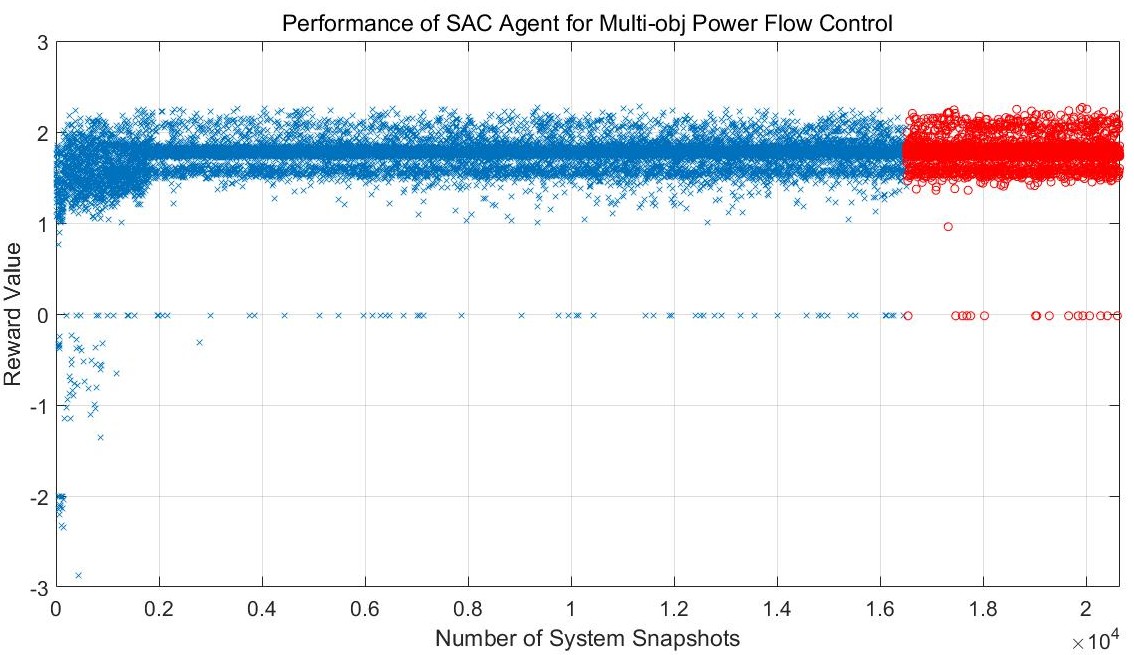}
  \caption{Performance of SAC agent (blue dots: training phase; red circle: testing phase)}
  \vspace{-10pt}
  \label{fig:SACperformance} 
\end{figure}

The presented SAC-based methodology for multi-objective power flow control was developed in Python 3.6 using Tensorflow 1.14, which was deployed in the control center of SGCC Jiangsu Electric Power Company in Nov., 2019. The city-level power grid (220 kV and above) is used to demonstrate the SAC agents' performance, which consists of 45 substations, 5 power plants (with 12 generators) and around 100 transmission lines, serving electricity to the city of Zhangjiagang, China. Massive historical operating snapshots (full topology node/breaker models for Jiangsu province with ~1500 nodes and ~420 generators, at an interval of 5 minutes) were obtained from their EMS system. The control objectives are to minimize transmission losses (at least 0.5\% reduction) without violating bus voltages ([0.97-1.07] pu) and line flow limits (100\% of MVA rating), namely, Eqs.~\ref{eq_1},~\ref{eq_8} and ~\ref{eq_9}. Voltage setpoints of the 12 generators in 5 power plants are adjusted by the SAC agent every 5 minutes. The performance of training and testing SAC agents using a time series of actual system snapshots is illustrated in Fig~\ref{fig:SACperformance} and Fig~\ref{fig:network losses}, where positive rewards indicate successfully resolving voltage and thermal violation issues, defined in  Eq.~\ref{eq_8}, Eq.~\ref{eq_9}. From 12/3/2019 to 1/13/2020, 7,249 operating snapshots were collected. Two additional epochs of the original data sets were created and randomly shuffled to create a training set (80\%) and a test set (20\%). For the first ~150 snapshots, the SAC agent struggles to find effective policies (with negative reward values); however, achieves satisfactory performance thereafter. Three training processes are simultaneously launched and updated twice a week to ensure control performance. For real-time application during this period, the developed method provides valid controls for 99.41\% of these cases. The average line loss reduction is 3.6412\% (compared to the line loss value before control actions). There are 1,019 snapshots with voltage violations, in which SAC agent solves 1,014 snapshots completely and effectively mitigates the remaining 5 snapshots. 

\begin{figure}
  \centering
  \includegraphics[width=13.7cm,height=6cm]{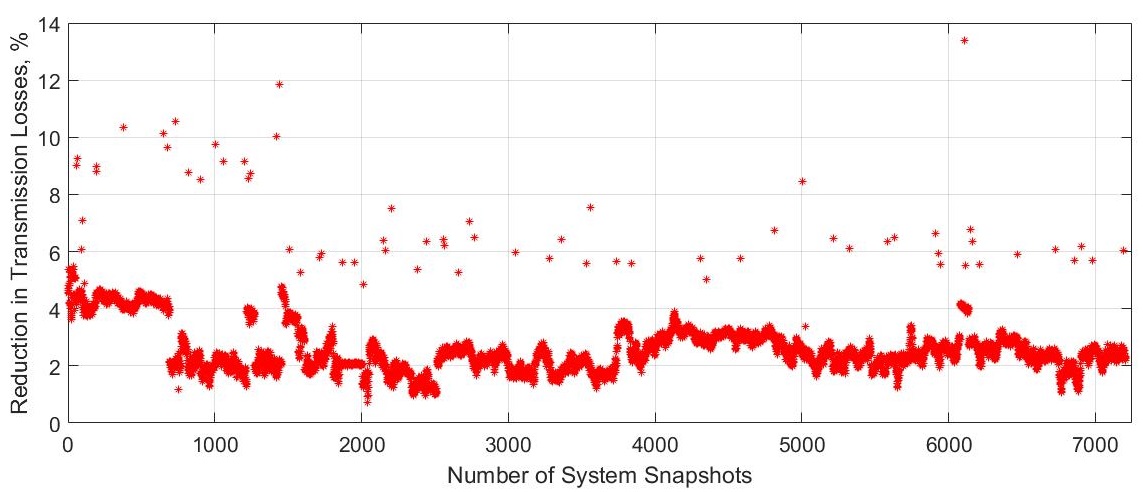}
  \caption{Reduction in transmission losses, City of Zhangjiagang}
  \vspace{-10pt}
  \label{fig:networklosses} 
\end{figure}

\section{Conclusions and Future Work}
In this work, a novel methodology of training effective SAC agents with periodic updating is presented for multi-objective power flow control. The detailed design and flowchart of the proposed methodology are provided for reducing transmission losses without violating voltage and line constraints. Massive numerical simulations conducted on a real power grid in real-time operational environment demonstrates the effectiveness and robustness of this approach. In future work, a multi-agent RL-based approach will be investigated for coordinating more types of corrective and preventive controls for use in larger power networks.

\section*{Broader Impact}
Power grid operation and control are facing new challenges caused by the increased penetration of intermittent energy and dynamic electricity demands. Many traditional planning and operational practices used by human operators may become ineffective under certain circumstances, because it is not possible to precisely predict and study all possible power grid operating conditions thoroughly for deriving corresponding controls, especially under severe disturbances that are never studied. In addition, the models across different time horizons used for grid planning and operational purposes are never perfect, resulting over-investment of grid infrastructure (if models are pessimistic) or security risks (if models are optimistic). 

The fast development and evolvement of AI technologies over the past few years and their successful stories reported in various control problems, provide a promising way to fill the above technology gaps facing today's power industry. In this work, we adopt state-of-the-art reinforcement learning algorithms to tackle multi-objective real-time power flow control problem and achieved preliminary success. During the development phase, several difficulties were encountered by the team, which may be typical for other industrial control problems, including (1) ensuring raw data and model quality to be used for training RL agents, (2) collecting a sufficient number of meaningful samples, (3) improving performance of RL agents via hyperparameter tuning with limited computing resources, (4) design of interfaces interacting with existing environment in a control center without interrupting existing system operation, (5) ensuring long-term control effectiveness and robustness, (6) coordination of AI-based controls with human operators. Besides the results presented in this paper, we are expanding the capability of the developed method toward larger areas in the provincial power network, with similar control performance observed. Moreover, a pilot project is underway in order to apply AI-based voltage controls in a closed-loop form, by switching capacitors and reactors for regulating voltage levels in a small zone using control commands derived from RL agents. 

We hope the findings of this work can benefit researchers and engineers in the area of promoting and applying AI techniques in the power community.   
\section*{Acknowledgments and Disclosure of Funding}
This work is funded by State Grid Corporation of China (SGCC) Jiangsu Electric Power Company through its Science and Technology Program.

\section*{References}
\small
[1]	E. Mohagheghi, M. Alramlawi, A. Gabash and P. Li, “A Survey of Real-Time Optimal Power Flow,” Energies, vol. 11, no. 11, 2018.

[2]	D. K. Molzahn, J. T. Holzer, B. C. Lesieutre, et al., “Implementation of a large-scale optimal power flow solver based on semidefinite programming,” IEEE Trans. Power Systems, vol. 28, no. 4, pp. 3987-3998, Nov. 2013.

[3]	R. Madani, M. Ashraphijuo, and J. Lavaei, “Promises of conic relaxation for contingency-constrained optimal power flow problem,” IEEE Trans. Power System, vol. 31, no. 2, pp. 199-211, Jan. 2015.

[4]	Y. Tang, K. Dvijotham, and S. Low, “Real time optimal power flow,” IEEE Trans. Smart Grid, vol. 8, no. 6, pp. 2963–2973, Nov. 2017.

[5]	DOE ARPA-E, “Grid Optimization Competition.” [Online]. Available: https://gocompetition.energy.gov/ 

[6]	R. Diao, Z. Wang, D. Shi, etc., “Autonomous voltage control for grid operation using deep reinforcement learning,” IEEE PES General Meeting, Atlanta, GA, USA, 2019.

[7]	B. Zhang, X. Lu, R. Diao, et. al. “Real-time Autonomous Line Flow Control Using Proximal Policy Optimization,” paper accepted for the IEEE PES General Meeting, Montreal, Canada, 2020.

[8]	Z. Yan and Y. Xu, “Data-driven load frequency control for stochastic power systems: A deep reinforcement learning method with continuous action search,” IEEE Trans. Power Syst., pp. 1–1, 2019.

[9]	Q. Huang, R. Huang, W. Hao, etc., “Adaptive power system emergency control using deep reinforcement learning,” IEEE Trans. Smart Grid, early access, 2019.

[10]	T. Haarnoja et al., “Soft actor-critic: off policy maximum entropy deep reinforcement learning with a stochastic actor,” in ICML, vol. 80, Stockholm Sweden, Jul. 2018, pp. 1861–187.

[11]	P. Kundur, Power System Stability and Control, McGraw-Hill, Inc, New York (1994)

[12]	T. Haarnoia, et al., “Soft actor-critic algorithms and applications,” arXiv preprint arXiv:1812.05905, 2018.

\section*{Appendix: }

\begin{table}[ht]
  \caption{Typical hyperparameters used for training SAC agents}
  \label{hyperparameter-table}
  \centering
  \begin{tabular}{lll}
    \toprule
    Name     & Description     & Value \\
    \midrule
    batch\_size & batch size  & 32-200     \\
    maxlen    & max number of iterations per episode & 10-40      \\
    Decay\_rate     & decay rate       & 0.15  \\
    replay\_size     & replay buffer size   & 100,000-1,000,000 \\
    gamma     & discount factor   & 0.99 \\
    epochs    & number of epochs & 1-5 \\
    polyak    & polyak   &  0.995 \\
    lr    & learning rate   &  0.0001-0.001 \\
    alpha    & initial value of alpha   &  0.001-0.2 \\
    entropy\_max    & max value of entropy   &  0.001 \\
    random\_seed    & random seed value   & 10-30 \\
    n\_hidden       & number of hidden layers  & 2-5 \\
    n\_neuros       & number of neuros in hidden layers & 64  \\
    n\_inputlayer   & dimension of input layer & 300-500\\
    n\_outputlayer  & dimension of output layer & 5\\
    \bottomrule
  \end{tabular}
\end{table}

\end{document}